\documentclass[12pt]{amsart}
\usepackage[square,compress,comma, numbers,sort]{natbib}
\usepackage[colorlinks=true, citecolor=blue, linkcolor=blue]{hyperref}
\usepackage{amsfonts}
\allowdisplaybreaks[4]
\usepackage{graphicx}
\usepackage{epstopdf}

\usepackage{amssymb}
\usepackage{amsmath}
\usepackage{color}
\usepackage{enumerate}
\newcommand{\abs}[1]{\left\lvert #1 \right\rvert}

\def\E#1{\mathbb{E}\left \{#1 \right\}}

\definecolor{c20}{rgb}{0.,0.7,0.}
\definecolor{c30}{rgb}{0.,0.,1.}
\definecolor{c40}{rgb}{1,0.1,0.7}
\definecolor{c50}{rgb}{1,0,0}
\definecolor{c60}{rgb}{1,0.9,0.1}
\definecolor{c70}{rgb}{0.50,1.00,0.00}

\def\Ehb#1{{\textcolor{c30}{#1}}}

\def\zE#1{{\textcolor{c30}{#1}}}

\def\N{\mathbb{N}}
\numberwithin{equation}{section}
\newtheorem{theo}{Theorem}[section]
\newtheorem{sat}[theo]{Proposition}
\newtheorem{de}[theo]{Definition}
\newtheorem{lem}{Lemma}[section]

\newtheorem{korr}[theo]{Corollary}
\newtheorem{remark}[theo]{Remark}
\newtheorem{remarks}[theo]{Remarks}
\numberwithin{equation}{section}

\newcommand{\prooftheo}[1]{ \textsc{Proof of Theorem} \ref{#1} }

\newcommand{\prooflem}[1]{\textsc{Proof of Lemma} \ref{#1}}

\newcommand{\pk}[1]{\mathbb{P} \left\{ #1 \right\} }

\newcommand{\QED}{\hfill $\Box$}
\newcommand{\COM}[1]{}
\def\IF{\infty}
\newcommand{\R}{\mathbb{R}}
\newcommand{\inr}{\in \R}

\newcommand{\BQN}{\begin{eqnarray}}
\newcommand{\EQN}{\end{eqnarray}}
\newcommand{\BQNY}{\begin{eqnarray*}}
\newcommand{\EQNY}{\end{eqnarray*}}
\def\ldot{, \ldots,}

\def\polhk#1{\setbox0=\hbox{#1}{\ooalign{\hidewidth
\lower1.5ex\hbox{`}\hidewidth\crcr\unhbox0}}}

\newcommand{\kb}[1]{\boldsymbol{#1}}
\newcommand{\vk}[1]{\kb{#1}}

\usepackage{mathtools}

\newcommand{\norm}[1]{\lVert #1 \rVert}
\newcommand{\nelem}[1]{{Lemma \ref{#1}}}

\newcommand{\netheo}[1]{{Theorem \ref{#1}}}

\def\vp{\varepsilon}

\def\IF{\infty}

\def\Cov{\mathrm{Cov}}

\date{}

\def\oo{(1+o(1))}

\def\LT{\left}
\def\RT{\right}

\def\Var{\text{Var}}

\newcommand{\limit}[1]{\lim_{#1 \to \infty}}

\newcommand{\BS}{\begin{sat}}
\newcommand{\ES}{\end{sat}}
\newcommand{\BT}{\begin{theo}}
\newcommand{\ET}{\end{theo}}
\newcommand{\BK}{\begin{korr}}
\newcommand{\EK}{\end{korr}}

\newcommand{\BD}{\begin{de}}
\newcommand{\ED}{\end{de}}
\newcommand{\BIT}{\begin{itemize}}
\newcommand{\EIT}{\end{itemize}}
\newcommand{\BDI}{\begin{description}}
\newcommand{\EDI}{\end{description}}
\newcommand{\BRM}{\begin{remarks}}
\newcommand{\ERM}{\end{remarks}}
\newcommand{\BEL}{\begin{lem}}
\newcommand{\EEL}{\end{lem}}

\def\LT{\left}
\def\RT{\right}

\def\Var{\text{Var}}

\def\Cov{\mathrm{Cov}}

\def\oo{(1+o(1))}

\def\LT{\left}
\def\RT{\right}

\def\Var{\text{Var}}

\def\oo{(1+o(1))}

\def\LT{\left}
\def\RT{\right}

\def\Var{\text{Var}}

\def\nj#1{ \mathbb{I}_u \left( #1 \right)}
\def\njk#1{ \mathbb{I}_0 \left( #1 \right)}
\def\njb#1{ \mathbb{I}\left( #1 \right)}
\def\MB{\mathcal{B}}

\def\Hal{\mathcal{H}_{\alpha}}

\def\Z{\mathbb{Z}}

  \def\td{\text{\rm d}}
\newcommand{\argmax}{{\rm argmax}}

\begin{document}
	
\title{Approximation of Sojourn Times of Gaussian Processes}
	
	\author{Krzysztof D\c{e}bicki}
	\address{Krzysztof D\c{e}bicki, Mathematical Institute, University of Wroc\l aw, pl. Grunwaldzki 2/4, 50-384 Wroc\l aw, Poland}
	\email{Krzysztof.Debicki@math.uni.wroc.pl}

	\author{Enkelejd  Hashorva}
	\address{Enkelejd Hashorva, Department of Actuarial Science,
		\\Faculty of Business and Economics\\
		University of Lausanne,\\
		UNIL-Dorigny, 1015 Lausanne, Switzerland
	}
	\email{Enkelejd.Hashorva@unil.ch}

	\author{Xiaofan Peng}
	\address{Xiaofan Peng, School of Mathematical Sciences, University of Electronic Science and Technology of China, Chengdu 610054, China}
	\email{xfpengnk@126.com}

\author{Zbigniew Michna}
	\address{Zbigniew Michna, Department of Mathematics, Wroc\l aw University of Economics, Poland}
	\email{Zbigniew.Michna@ue.wroc.pl}

	\bigskip
	
	\date{\today}
	\maketitle

{\bf Abstract}:  We investigate the tail asymptotic behavior of the sojourn time for a large class of centered Gaussian processes $X$,
in both continuous- and discrete-time framework.
All results obtained here are new for the discrete-time case.
In the continuous-time case, we complement the investigations of \cite{berman1987extreme,MR803245} for non-stationary $X$.
A by-product of our investigation is a new representation
of Pickands constant which is important for Monte-Carlo simulations and yields a sharp lower bound for Pickands constant.



\section{Introduction}
Let $X(t),t\inr $ be a centered Gaussian process with variance function $\sigma^2$,  correlation function $\rho$ and continuous trajectories.
By
\BQNY\label{defi-sojo}
L_u[a,b] :=\int_{a}^b \nj{X(t)}\td t
\EQNY
we define
the sojourn time spent above a fixed level $u$ by the process $X$ on the interval $[a,b]$,
where $\nj{x}:= \njb{x> u}$.

In a series of papers culminating in \cite{Berman92},
S. Berman derived results on the tail asymptotic behaviour
of  $L_u[a,b]$, as $u\to \IF$.
{\it The sojourn time approach} to tackle this problem consists in finding explicitly an appropriate scaling function $v(u)$ such that
for some function $C(x)>0, x\ge 0$
\BQN \label{CX}
\pk{ L_u[a,b] > x/v(u)} &\sim& C(x) \pk{ L_u[a,b] > 0}= C(x) \pk{\sup_{t\in [a,b]} X(t)> u} ,\quad u\to\IF
\EQN
for  any $x\ge 0$ a continuity point of $C(\cdot)$.
In our notation $\sim$ stands for asymptotic equivalence of two
functions as the argument tends to 0 or $\IF$.
Additional inside of this approach is the explicit calculation of the
exact asymptotics of
$\pk{\sup_{t\in [a,b]} X(t)> u}$ as $u\to \IF$.
For example, as shown in several works of Berman and Pickands (see e.g.,
\cite{Berman82,Berman92, PickandsB}) for $X$ a centered stationary Gaussian process the asymptotic tail behaviour
of $v(u)L_u([a,b])$ and that of $\sup_{t\in [a,b]} X(t)$ can be studied under appropriate assumptions on the correlation function $\rho$.
Pickands' assumption for $X$ stationary with unit variance function  reads
\BQN\label{Pic}
1- \rho(t) \sim \abs{t}^\alpha, \ t\to 0\quad \textrm{and}\quad \rho(t) < 1,\ \forall\, t\not=0,
\EQN
where $\alpha \in (0,2]$. Under \eqref{Pic} in view of \cite{PickandsB} (see also \cite{Pit20}) taking the scaling function $ v(u)= u^{2/\alpha}$
we have (consider for simplicity $[a,b]=[0,T],T>0$)
\BQN\label{eq:etaA}
 \pk{\sup_{t\in [0,T]} X(t)> u} \sim T \Hal v(u) \pk{X(0)> u}, \quad u\to \IF,
\EQN
where $\Hal $ is the Pickands constant given by
$$\Hal =\lim_{S\to\IF}S^{-1}\Hal([0,S])\in (0,\IF),$$
with
\BQN\label{def-H-al-e}
 \Hal([0,S])=\E{\sup_{t\in  [0,S]}e^{W_\alpha(t)}} ,
\quad W_\alpha(t):=\sqrt{2}B_{\alpha}(t)-\abs{t}^{\alpha}
\EQN
and $B_{\alpha}$  is a standard fractional Brownian motion (fBm) with Hurst index $\alpha/2 \in (0,1]$.
 A refinement  of \eqref{eq:etaA} is given in \cite{Berman92}[Theorem 3.3.1]. Namely, \eqref{CX} holds with
 $$C(x)=\widetilde \MB_{\alpha} (x)/\Hal$$
 for any $x>0$ a continuity point of $\widetilde \MB_\alpha(\cdot) $. Here
$ \widetilde \MB_{\alpha}(x) = \int_{x}^{\IF}\frac{1}{y}\, \td\mathcal{G}_{\alpha}(y)\in (0,\IF), $ with
\BQN\label{defi-g-ale-x}
\mathcal{G}_{\alpha}(x) = \pk{ \int_{ \R } \njk{W_{\alpha}(s)  +  \mathcal{E}} \td s \leq x },
\EQN
where $ \mathcal{E}$  is a unit exponential random variable independent of $W_\alpha$.
Furthermore,  as shown in \cite{Berman92}[Theorem 10.5.1]
 \BQN \label{H0B}
\widetilde \MB_\alpha (0)= \lim_{x\downarrow 0} \int_x^\IF \frac{1}{y}\, \td G_{\alpha}(y) = \mathcal{H}_\alpha\,.
\EQN

We note that the only known values of Pickands constants are
$\mathcal{H}_1=1$ and $\mathcal{H}_2=\frac{1}{\sqrt{\pi}}$
and both \eqref{def-H-al-e} and \eqref{H0B} are not
tractable for simulations.
In Theorem \ref{BermanPic} we present an interesting formula for $\mathcal{H}_\alpha$, which is a consequence of Berman's theory on extremes of random processes.
We believe that this new formula is of particular interest for simulations, since it is given as an expectation, see \cite{MR2458013,DiekerY,DM, SBK,KrzysEH} for alternative formulas. Another advantage of this new formula is that it implies the uniformly (with respect to $\alpha$) sharpest lower bound for the Pickands
constant  available in the literature so far.
Next, let $\Gamma(\cdot)$ stands for Euler Gamma function.

\BT \label{BermanPic} For $\alpha \in (0,2]$  we have
\BQN\label{BermanRep}
\mathcal{H}_\alpha  = \E{ \frac{1}{
		\int_{ \R} \njk{W_\alpha(s) + \mathcal{E} } ds }}  
\ge\frac{\Gamma(1/\alpha)}{4\Gamma(2/\alpha)}.
\EQN
\ET

\begin{figure}[!ttt]
\includegraphics[width=11cm]{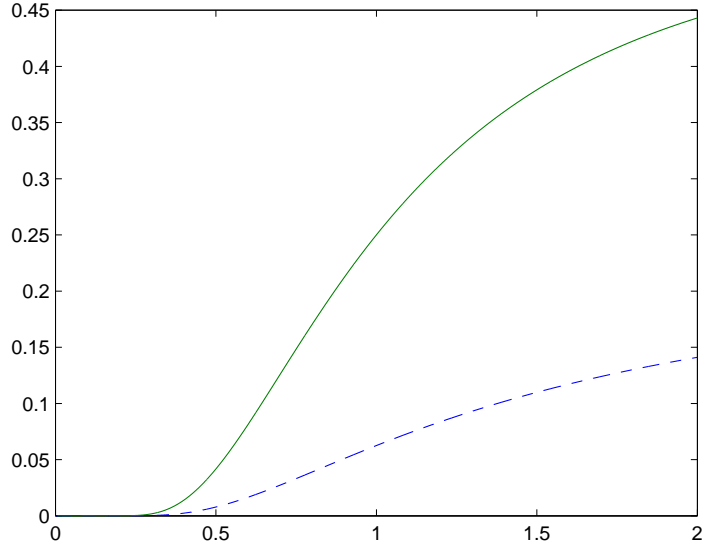}
	\caption{Solid line: lower bound (\ref{BermanRep}),
dashed line: lower bound (\ref{dmr})}\label{Fig}
\end{figure}
Interestingly, the same lower bound for $\mathcal{H}_\alpha$
as derived in Theorem \ref{BermanPic}
was
obtained heuristically in \cite{Aldous}[J20a,J20b].
The above finding uniformly improves the
result of
 \cite{mi:17} (see also \cite{MR1993949}):
\BQN\label{dmr}
\mathcal{H}_\alpha \ge  \frac{4^{-\frac{1}{\alpha}-1}}{\Gamma(1/\alpha+1)};
\EQN
see Fig. \ref{Fig}.
We refer to \cite{Harper2} for the proof that
$\mathcal{H}_\alpha \ge  \frac{ (1.1527)^{1/\alpha}}{ \Gamma(1/\alpha)}$
for $\alpha$ sufficiently close to $0$,
which subverted an opened for long time hypothesis that
$\mathcal{H}_\alpha=\frac{1}{\Gamma(1/\alpha)}$.
Other estimates for Pickands constants can be found in e.g., \cite{sh:96} and
\cite{MR2458013}.

The main interest of this contribution is the investigation of the tail asymptotics of $L_u[a,b]$ and its discrete counterpart.\\
Our method here is completely different from that of Berman.
Namely, in this paper we developed the uniform double-sum method for the sojourn time functional.
Interestingly, this approach leads to a new representation of Berman's constants $\widetilde \MB_{\alpha}(x)$; see Section \ref{sect-main-re}
where the asymptotics for the tail distribution of sojourns of locally-stationary Gaussian processes was derived and
compared with the classical results of Berman.

Our main findings in this paper can be summarized as follows: for both locally-stationary Gaussian processes
and general non-stationary Gaussian processes with variance maximal at some unique point, we show that \eqref{CX} holds for almost all $x$ and moreover, we calculate explicitly $C(x)$ and give the appropriate scaling function $v$.  Our results are new for non-stationary Gaussian processes, and agree with those of Berman for the locally stationary ones. In particular, all results are new for the discrete setup introduced in the next section.

Brief organisation of the rest of the paper:
{In Section 2 we derive the tail asymptotics of sojourn time for locally stationary Gaussian processes. Corresponding results for general non-stationary Gaussian processes are then presented in Section 3. All the proofs are displayed  in Section 4 whereas few technical results are included in Section  5.}

\section{Sojourns of Locally Stationary Gaussian Processes} \label{sect-main-re}
In this section we analyze sojourns for the class
of {\it locally stationary} Gaussian processes, introduced by Berman in \cite{Berman92}, see also \cite{MR1342094,MR1062062,MR3679987, ChanLai,Marek,MR3638374}.
Specifically, let $X(t),t\in [0,T]$ be a centered Gaussian process with unit variance and correlation function $\rho$ satisfying
 \BQN\label{locall-station}
\lim_{\epsilon\to0}\sup_{t,t+s\in[0,T], |s|<\epsilon} \abs{ \frac{1-\rho(t,t+s)} {K(\abs{s})} - H(t) }=0,
\EQN
where $H$ is a continuous positive function on $[0,T]$ and $K$ is a regularly varying function at 0 with  index $\alpha\in(0,2]$. In the following let $v$ be the asymptotically unique function (which exists, see \cite{Berman92}) such that $\limit{u} v(u)=\IF$  and
 \BQN\label{defini-vu-loc}
\lim_{u\to\IF} u^2 K(1/v(u))=1.
\EQN
We shall investigate the tail asymptotics of $L_u^*[0,T]:=v(u) L_u[0,T]$.
 Given some $\eta>0$ we define the discrete counterpart of $L_u^*[a,b]$ as
$$ L_{\eta, u}^*[a,b]:=  v(u)\int_{a}^{b}  \nj{X(t)} \mu_{\eta_u}(\td t)
= \eta \sum_{  t\in(\eta_u \Z) \cap [a,b]  } \nj{X(t)} ,$$
where $\eta_u= \eta/v(u)$ and $\mu_c(dt)/c$ denotes the counting measure on $c\Z$. In the sequel we interpret $0 \Z$ as $\R$ and $\mu_0$ as the Lebesgue measure on $\R$. Since $ \mu_{c}$ converges to the Lebesgue measure $\mu_0$ on $\R$ as $c\to 0$,  with this convention we set
$$L_{0,u}^*[a,b]:=L_u^*[a,b]= v(u) L_u[a,b] .$$

In order to state our first result, define  for any $\lambda>0, \eta \ge 0, x\in[0,\mu_\eta([0,S]))$
\BQN\label{Baleta}
 \MB_{\alpha,\lambda}^\eta(S,x):= \int_{\R} \pk{ \int_{0}^{S} \njk{W_\alpha(\lambda^{1/\alpha} s) +z } \mu_\eta(\td s) >x } e^{-z} \td z,
 \EQN
where $W_\alpha$ is defined in \eqref{def-H-al-e}.
Further, for any $x\ge 0$ set
\BQN\label{def-B-H-x}
\MB_{\alpha}^{\eta,H} (x) := \lim_{S\to\IF} \frac{\int_0^T \mathcal{B}_{\alpha,H(t)}^\eta (S,x)\td t}{S}, \quad \MB_{\alpha}^\eta (x) :=\lim_{S\to\IF} S^{-1}\MB_{\alpha,1}^\eta (S,x).
\EQN
Hereafter, when we mention that $x$ is a continuity point for some function $f$
 {we also assume  that $f(x)>0$}.

Next we state our first result. The case $\eta>0$ is new, whereas for the case $\eta=0$ we retrieve the result of Berman,
however the asymptotic constant (pre-factor) is given in a different form
than in the original Berman's result, see e.g. \cite{Berman82}, which is due
to a different technique applied here.

\BT\label{theo-ls-tail}
Let $X(t),t\in [0,T]$ be a centered, sample path continuous Gaussian process with unit variance and correlation function satisfying assumption \eqref{locall-station}. If further $\rho(s,t)<1$ for all $s,t\in[0,T],s\neq t$,
then for any $x>0$ a continuity point of $\MB_{\alpha}^{\eta,H}(\cdot)$ and for $x=0$ we have
\BQN\label{concl-sta-tail}
\pk{L_{\eta,u}^*[0,T]>x} \sim \MB_{\alpha}^{\eta,H}(x) v(u) \pk{X(0)> u}, \quad u\to\IF,
\EQN
where $v(u)$ is given in \eqref{defini-vu-loc} and $\MB_{\alpha}^{\eta,H}(\cdot)$
defined in \eqref{def-B-H-x}  is positive and finite for any $x,\eta\ge 0$.
\ET

\begin{remark} \label{remark-p-b}
i) If $X(t),t\in [0,T]$ is a centered, stationary, sample path continuous Gaussian process with unit variance function
and its correlation function $\rho$ satisfies Pickands condition \eqref{Pic}, then
$X$ is locally stationary with function  $H(t)\equiv1,t\in [0,T]$.  For such $H$ we have that
$\MB_{\alpha}^{\eta,H}(x) = T \MB_{\alpha}^{\eta}(x)$.\\
ii) For $\eta=0$,  by \cite{Berman92}[Theorem 3.3.1] and \eqref{concl-sta-tail}
we have
\BQNY\label{densS} \MB_\alpha^0(x)=\widetilde \MB_{\alpha}(x)
\EQNY
for all continuity points of $\widetilde \MB_{\alpha}(\cdot)$
(since both $ \MB_\alpha^0(\cdot)$ and $\widetilde\MB_{\alpha}(\cdot)$ are monotone non-increasing).
\end{remark}

\section{Sojourns of Non-Stationary Gaussian Processes}
In this section we analyze sojourns of non-stationary centered Gaussian processes.
Suppose that $X(t),t\in [-T,T]$ is a centered Gaussian process with continuous sample paths.
Tractable assumptions
on both variance $\sigma^2(t)=Var(X(t))$ and correlation function $\rho(s,t)$,
adopted
from a vast literature on the asymptotic analysis of
supremum of non-stationary Gaussian processes,
see e.g., \cite{Pit72,Pit96,MR823082,MR803245,MR985382,berman1987extreme,dbicki2015gaussian,Pit20,MR3679987,high2},
are as follows:
	\begin{itemize}
		\item[{\bf A0}:] For some $T>0$
		$$t_0=\argmax_{t\in [-T,T]} \sigma(t)$$
		is unique. For notational simplicity we assume further that
	$t_0=0$ and  $\sigma(t_0)=1.$
		\item[{\bf A1}:] For some $\alpha\in (0,2]$ we have
		\BQNY\label{assump-corre}
		1- \rho(s,t)\sim \abs{t-s}^\alpha,\qquad s,t\to t_0.
		\EQNY	

		\item[{\bf A2}:] For some positive constants $b,\beta$
		\BQNY\label{assump-vari}
		1- \sigma(t)\sim b\abs{t}^\beta 
		,\qquad t\to t_0.
		\EQNY
	\end{itemize}
	
Under the assumptions {\bf A0-A1}, if further
\BQNY
\lim_{s,t \to 0,s\neq t} \frac{ \abs{\sigma(s)- \sigma(t)}}{ \E{ ( X(s)- X(t) ) ^2}}=0
\EQNY
in view of \cite{berman1987extreme}[Theorem 6.1] we have
\BQNY
\limit{u} \frac{\int_0^x y \td \pk{ L_{u }^*[-T,T] \le y}}{\E{L_{ u}^*[-T,T]}} = \mathcal{G}_\alpha (x)
\EQNY
for any continuity point $x$ of $\mathcal{G}_\alpha $ defined in \eqref{defi-g-ale-x}.\\
See also \cite{MR803245} for another result shown under {\bf A0, A2} assuming further that
\BQNY
\lim_{t \to 0} \frac{ \E{ ( X(t)-X(0) )^2}}{ 1-  \sigma(t)}=0.
\EQNY

Under the assumptions \textbf{A0-A2} we shall derive the tail asymptotics of $L_{\eta, u}^*[-T,T]$, where we chose the scaling function $v(u)$
as follows
 \BQN\label{defi-v}
v(u)=u^{2/\min(\alpha,\beta)}.
\EQN
As in the case of Piterbarg's result for $\sup_{t \in [-T,T]} X(t)$ (see \cite{Pit20}), if $\alpha=\beta$ in the asymptotic results a new constant  $\mathcal{P}_\alpha^{b,\eta} $ appears, which is defined for any $b>0$ by
\BQN\label{defi-B-alp-b}
\mathcal{P}^{b,\eta}_\alpha(x) := \int_{\R} \pk{\int_{\R }
	\njk{ W_{\alpha}(s)-b\abs{s}^\alpha+ z }\mu_\eta(ds)>x } e^{-z}\td z.
\EQN

Additionally, for $\eta>0$
we set
\BQN\label{conT1}
\mathcal{T}^{b,\eta}_\beta(x)
=
\left\{
\begin{array}{ll}
1 & \mbox{if $x\in [0,\eta) $}\\[0.5cm]
e^{-b(k\eta)^\beta} & \mbox{if $x\in[(2k-1)\eta,(2k+1)\eta), k\in \N $}
\end{array}
\right.
\EQN
and for $\eta=0$, $x\geq0$ let
$\mathcal{T}^{b,0}_\beta (x)=e^{-b(\frac{x}{2})^\beta}.$

We present next the main result of this section.

\BT \label{MainTheo} Let $X$ be a centered Gaussian process satisfying {\bf A0-A2},
$\eta\ge0$
and  $v(u)=u^{2/\min(\alpha,\beta)}$.
\\
i) If $\alpha<\beta$, then for any $x>0$ a continuity point of $\MB_{\alpha}^\eta (\cdot)$ 
and $x=0$  
\BQN\label{con-alp-s-bet-tail}
\pk{L_{\eta,u}^*[-T,T] >x}\ \sim \ 2b^{-1/\beta}\Gamma(1/\beta+1)\MB_{\alpha}^\eta (x) u^{2/\alpha-2/\beta}\pk{X(0)> u}  , \quad u\to \IF.
\EQN
ii)  If $ \alpha=\beta$, then for any $x>0$ a continuity point of  $\mathcal{P}_\alpha^{b,\eta}(\cdot)$ and $x=0$
\BQN \label{LCohen}
\pk{L_{\eta,u}^*[-T,T]>x} \ \sim\  \mathcal{P}^{b,\eta }_\alpha(x)\pk{X(0)> u} ,
\quad u\to \IF.
\EQN
iii)  If $ \beta < \alpha $,
then for any $x>0$ a continuity point of  $\mathcal{T}^{b,\eta}_\beta(\cdot) $  and $x=0$
\BQN \label{con-T-beta}
\pk{L_{\eta,u}^*[-T,T]>x} \ \sim\ \mathcal{T}^{b,\eta}_\beta(x)  \pk{X(0)> u} ,
\quad u\to \IF.
\EQN
\ET
\begin{figure}[!ttt]
\includegraphics[width=11cm]{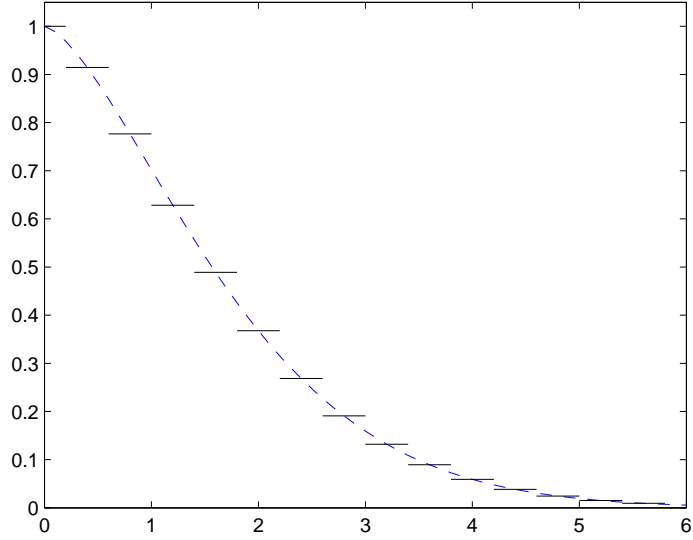}
	\caption{Graph of  $\mathcal{T}_{1.5}^{1,\eta}(\cdot)$. Solid line: $\eta=0.2$, dashed line:  $\eta=0$. }\label{Fig2}
\end{figure}

\begin{remark}\label{remar-1}
	i) If $\eta=x=0$, then
$\mathcal{P}^{b,0}_\alpha(0)=\mathcal{P}_{\alpha}^b$.
 Indeed, for any $b>0$ we have
	\BQNY
	\mathcal{P}^{b,0}_\alpha(0) &=& \int_{\R} \pk{\int_{s \in \R }
		\njk{ W_{\alpha}(s)-b\abs{s}^\alpha+ z }\td s>0 } e^{-z}\td z\\
	&=& \int_{\R} \pk{ \sup_{s\in\R} \Bigl(W_{\alpha}(s)-b\abs{s}^\alpha+ z\Bigr) >0 } e^{-z}\td z\\
	&=& \E{\sup_{s\in \R}e^ {W_{\alpha}(s)-b\abs{s}^\alpha } }=: \mathcal{P}_{\alpha}^{b}.
	\EQNY
In the literature, $\mathcal{P}_{\alpha}^{b}$ is referred to as Piterbarg constant, see \cite{Pit96, MR1814702, KEP16,PBConst} for related constants and basic properties. \\
ii) 
If $t_0\in \{-T,T\}$ in {\bf A0}, then \netheo{MainTheo} still holds  subject to appropriate change of the constants
in the asymptotics.
Specifically, if $\alpha<\beta$, \eqref{con-alp-s-bet-tail} holds with the constant
$2$ removed from the expression. If $\alpha=\beta$, then $\mathcal{P}^{b,\eta}_\alpha(x)$ in \eqref{LCohen}
has to be changed to
\BQNY\label{defi-B-alp-b2}
\int_{\R} \pk{\int_{0}^\infty
	\njk{ W_{\alpha}(s)-b\abs{s}^{\alpha}+ z } \mu_\eta(\td s)>x } e^{-z}\td z.
\EQNY
If $\alpha>\beta$, then in \eqref{con-T-beta}  $\mathcal{T}^{b,\eta }_\beta(x)$ has to be substituted by
$e^{-bx^\beta},x\geq0$ for $\eta=0$, and
$e^{-b(k\eta)^\beta}$ if $x\in[k\eta,(k+1)\eta), k\in \{0\}\cup\N $ for $\eta>0$.

\end{remark}

\def\vv{\nu}
\def\LL{\rho}

\section{Proofs}
Below $\lfloor x \rfloor$ stands for the integer part of $x$
{and $\lceil x\rceil$ is the smallest integer not less than $x$}. Further $\Psi$ is the survival function of an $N(0,1)$ random variable.

\def\Bale{\mathcal{B}_{\alpha,\lambda}^\eta}
\def\Seu{S_{\eta,u}}
\def\Se{S_\eta}

\prooftheo{BermanPic}
Since $W_\alpha$ has almost surely continuous trajectories with $W_\alpha(0)= 0$ and $\mathcal{E}>0$ almost surely,  then $I_\alpha=\int_{ \R} \njk{W_\alpha(s) + \mathcal{E} }ds
>0$ almost surely.  Consequently, by the definition of Pickands constant in \eqref{H0B}  and the monotone convergence theorem we obtain
$$ \mathcal{H}_\alpha=  \lim_{x\downarrow 0} \int_x^\IF \frac 1 y \, \td \mathcal{G}_\alpha(y)
= \E{ \frac 1 {{I}_\alpha}} \in (0,\IF).  $$
Hence by Jensen's inequality  we have
\BQNY \E{ \frac 1 {{I}_\alpha}}  &\ge & \frac 1 {\E{{I}_\alpha}}.
\EQNY
Further, we have
\BQN\nonumber
\E{ I_\alpha} &=& \int_{\R} \pk{ W_\alpha(t)+ \mathcal{E}>0  } \, \td  t \\\nonumber
&=& 2\int_0^\IF  \pk{ W_\alpha(t)+ \mathcal{E}>0  } \, \td  t \\\nonumber
&=& 2\int_0^\IF  \pk{ \sqrt{2 t^{\alpha}} B_\alpha(1)- t^\alpha + \mathcal{E}>0  } \, \td  t \\
&=& 4\int_0^\IF \pk{ B_\alpha(1)> \sqrt{\frac{t^\alpha}{2}}} \, \td  t \label{l4-proof}\\\nonumber
&=& 4\int_0^\IF \pk{ \max(0,B_\alpha(1))> \sqrt{\frac{t^\alpha}{2}}} \, \td  t\\\nonumber
&=& 4\E{ 2^{1/\alpha}({\max(0,B_\alpha(1))})^{2/\alpha}}\\\nonumber
&=& \frac{4^{1/\alpha+1/2}}{\sqrt{\pi}}\Gamma(1/\alpha+1/2),
\EQN
where in (\ref{l4-proof}) we used Lemma \ref{Tilt} from Appendix. Thus the proof is complete.
\QED

\prooftheo{theo-ls-tail}
Let $S>1$ be a positive constant. Define $S_0=S, \Se= \eta \lfloor S \rfloor $ for $\eta>0$ and set further
$$\Delta_k=[k \Seu,(k+1) \Seu],\quad  k=0 ,\ldots,N_u,
$$
where $\Seu=S_\eta/v(u)$ and $N_u=\lfloor T /\Seu \rfloor$.  We have for all $u$ positive and $x\ge 0$
\BQN\label{low-upp-I1}
 I_1(u) \leq \pk{ L_{\eta,u}^*[0,T] >x} \leq  I_2(u),
\EQN
 where
\BQNY
&& I_1(u) = \sum_{k=0}^{N_u-1}\pk{ L_{\eta,u}^*\Delta_k>x  } - \sum_{0\leq i<k\leq N_u-1}
q_{i,k}(u),\\
&&  I_2(u) = \sum_{k=0}^{N_u}\pk{ L_{\eta,u}^*\Delta_k>x } + \sum_{0\leq i<k\leq N_u}
q_{i,k}(u),
\EQNY
with
\BQNY \label{qik}
q_{i,k}(u)=\pk{ \sup_{t\in\Delta_i} X(t)>u, \sup_{t\in\Delta_k} X(t)>u }.
\EQNY
We first show that, as $u\to\IF$ and then $S\to\IF$, the first sum in $I_1$ is asymptotically equivalent to $v(u)\Psi(u)$ and the double sum is negligible with respect to the former one.\\
For any $x\geq0$ and $t\in[0,T]$, put
$$F_u(t,x)= \Psi^{-1}(u)\pk{ v(u)\int_0^{\Seu}  \nj{X(t+s)} \mu_{\eta_u}(\td s)>x}.$$
According to \eqref{locall-station}, we choose $\vp$ small enough such that
\BQN\label{defi-vp-con1}
1-\rho(s,t)\leq 2\overline h K(\abs{t-s})\leq1,\quad \forall s,t\in[0,T],\abs{t-s}\leq\vp
\EQN
with $\overline h=\max_{t\in[0,T]}H(t)$. Let $Y$ be a centered stationary Gaussian process with continuous trajectories, unit variance function and
covariance function satisfying
$$1-\Cov(Y(t),Y(t+s))\sim 4\overline hK(\abs{s}),\quad  s\to0.$$
The existence of such a Gaussian process is guaranteed by the Assertion in \cite{VP1}[p.265] and follows from \cite{MR906871,MR0045966}.
 Consequently,  by Slepian lemma and \cite{KEP16}[Lemma 5.1] for any $\eta\geq0$ and sufficiently large $u$
\BQN\nonumber
\sup_{t\in[0,T]} F_u(t,x) &\leq& \Psi^{-1}(u) \sup_{t\in[0,T]}\pk{ \sup_{s\in[0,\Seu]}X(t+s)>u }\\\nonumber
&\leq& \Psi^{-1}(u)\pk{ \sup_{s\in[0,\Seu]}Y(s)>u }\\\label{uni-bou-F}
&\leq& 2\lceil(4{\overline h})^{1/\alpha}\Se\rceil \Hal([0,1]),
\EQN
where $\Hal(\cdot)$ is defined in \eqref{def-H-al-e}.
Therefore,
\BQNY
&&\abs{ \frac 1 {\Se }\int_0^T   \frac 1 { \Psi(u)}\pk{L_{\eta,u}^*[t,t+  \Seu]>x} \mu_{\Seu}(\td t)
- \frac 1 { \Psi(u)v (u)}\sum_{k=0}^{N_u-1}\pk{ L_{\eta,u}^*\Delta_k>x  } }\\
&\leq& \frac1{v(u)}\sup_{t\in[0,T]} F_u(t,x)\rightarrow 0
\EQNY
as $u\to\IF$ implying
\BQNY
\frac 1 { \Psi(u)v (u)}\sum_{k=0}^{N_u-1}\pk{ L_{\eta,u}^*\Delta_k>x  }
\sim 
\frac 1 {\Se } \int_0^T  F_u(t,x)  \mu_{\Seu}(\td t).
\EQNY
 Let $x_0\in(0,\mu_\eta([0,S_\eta]))$ be a continuity point of $\int_0^T\mathcal{B}_{\alpha,H(t)}^\eta (\Se,x)dt$, then
   $$\lim_{\vp_0\to0} \int_0^T[\mathcal{B}_{\alpha,H(t)}^\eta (\Se,x_0+\vp_0) -\mathcal{B}_{\alpha,H(t)}^\eta (\Se,x_0-\vp_0)]\td t=0.$$
   Application of the dominated convergence theorem with Lemma \ref{upp-bou-BaeDx} in Appendix yields
  		$$\int_0^T [ \mathcal{B}_{\alpha,H(t)}^\eta (\Se,x_0+)-\mathcal{B}_{\alpha,H(t)}^\eta (\Se,x_0-)]\td t=0.$$
  Since $\mathcal{B}_{\alpha,H(t)}^\eta (\Se,x)$ is monotone in $x$ for each $t\in[0,T]$, it follows that $x_0$ is a
  continuity point for any $t\in[0,T]\setminus B$, where $B$ is some subset of $[0,T]$ with Lebesgue measure $0$.
Next, by \nelem{lemm-soj-sta-sma}-i), for any $t_u$ such that $\limit{u} t_u=t_0\in [0,T]\setminus B$
\BQN\label{limit-F}
\lim_{u\to\IF} F_u(t_u,x_0) =  \mathcal{B}_{\alpha,H(t_0)}^\eta (\Se,x_0).
\EQN
  By \eqref{uni-bou-F}
for sufficiently large $u_0$,  $F_u(\cdot,x_0),u\geq u_0$ is uniformly bounded on $[0,T]$. Consequently, \cite{Htilt}[Lemma 9.3] implies
\BQN
\limit{u} \frac 1 { \Psi(u)v (u)}\sum_{k=0}^{N_u-1}\pk{ L_{\eta,u}^*\Delta_{k}>x_0  }  &=&
\limit{u}  \frac 1 {\Se }\int_0^T  F_u(t,x_0)  \mu_{\Seu}(\td t)\notag \\\label{wea-con-loc}
&=&		\frac{1}{\Se} \int_0^T \mathcal{B}_{\alpha,H(t)}^\eta (\Se,x_0)\td t.
\EQN
Further, by \nelem{lemm-soj-sta-sma}-i), \eqref{limit-F} is also valid for $x_0=0$. Therefore, \eqref{wea-con-loc} holds for $x_0=0$ and for any $x_0\in(0,\mu_\eta([0,S_\eta]))$ a continuity point of $\int_0^T\mathcal{B}_{\alpha,H(t)}^\eta (\Se,x)dt$.

Define
$$A_\vp=\{(i,k):1\leq i+1<k\leq N_u-1, k+1-i\leq \lfloor\vp/\Seu\rfloor\},$$
$$B_\vp=\{(i,k):1\leq i+1<k\leq N_u-1, k+1-i> \lfloor\vp/\Seu\rfloor\}.$$
Then
\BQN\label{dou-sum-ls}
\sum_{0\leq i<k\leq N_u-1} q_{i,k}(u) \leq
\sum_{0\leq i\leq N_u-1} q_{i,i+1}(u) + \sum_{(i,k)\in A_\vp} q_{i,k}(u) +
\sum_{(i,k)\in B_\vp} q_{i,k}(u).
\EQN
Since $\rho(s,t)<1$ for all $s,t\in[0,T],s\neq t$, with a similar argument as used in the proof of Theorem 4 in \cite{mi:17}
\BQN\label{negi-far-ls}
\limsup_{u\to\IF}\frac{1}{v(u)\Psi(u)}\sum_{(i,k)\in B_\vp} q_{i,k}(u)=0
\EQN
holds. In view of \nelem{lemm-negli-ds}, for large enough $u$
\BQNY
\sum_{(i,k)\in A_\vp} q_{i,k}(u) \leq 2\lceil (16\overline{h})^{1/\alpha}\rceil^2 \lceil \Se\rceil^2 H_{\alpha}^2([0,1])\Psi(u)
 N_u\sum_{k=1}^\IF  \exp(-\frac{1}{16}\underline{h}\abs{k\Se}^{\alpha/2}).
\EQNY
Since $\lim_{S\to\IF}e^S\sum_{k=1}^\IF e^{-Sk^{\alpha/2}} < 2 $, then for sufficiently large $S$
\BQN\label{negi-notfar-ls}
\limsup_{u\to\IF}\frac{1}{v(u)\Psi(u)}\sum_{(i,k)\in A_\vp} q_{i,k}(u)\leq  4\lceil (16\overline{h})^{1/\alpha}\rceil^2 H_{\alpha}^2([0,1])T \Se
\exp(-\frac{1}{16}\underline{h}\Se^{\alpha/2}).
\EQN
Further, choosing  large $S$ such that $\Se>1$, for large $u$ and each $i< N_u$ we have
\BQNY
&& \frac{q_{i,i+1}(u)}{\Psi(u)}\\
&&\ \leq \frac1{\Psi(u)}\pk{ \sup_{t\in[(i+1)\Se,(i+1)\Se+\sqrt{\Se}]/v(u)}X(t)>u  }\\
 &&\quad\quad +\ \frac1{\Psi(u)}\pk{\sup_{t\in\Delta_{i}}X(t)>u, \sup_{t\in[(i+1)\Se+\sqrt{\Se},(i+2)\Se]/v(u)}X(t)>u}\\
 &&\leq 2\lceil(4{\overline h})^{1/\alpha}\sqrt{\Se}\rceil H_{\alpha}([0,1])+ 2\lceil(16\overline{h})^{1/\alpha}\rceil^2 \lceil \Se\rceil \lceil \Se-\sqrt{\Se}\rceil H_{\alpha}^2([0,1]) \exp(-\frac{1}{16}\underline{h}\abs{\Se}^{\alpha/4}),
\EQNY
where in the last inequality we have used \eqref{uni-bou-F} and \nelem{lemm-negli-ds}. Therefore,
\BQN\nonumber
&&\limsup_{u\to\IF}\frac{1}{v(u)\Psi(u)}\sum_{i=0}^{N_u-1} q_{i,i+1}(u)\\\label{neg-near-part-ls}
&&\quad \leq  \frac{2TH_{\alpha}([0,1])}{\Se} \Big( \lceil(4{\overline h})^{1/\alpha}\sqrt{\Se}\rceil
+ \lceil(16\overline{h})^{1/\alpha}\rceil^2 \lceil \Se\rceil^2 H_{\alpha}([0,1]) \exp(-\frac{1}{16}\underline{h}\abs{\Se}^{\alpha/4}) \Big).
\EQN
Substituting \eqref{negi-far-ls}-\eqref{neg-near-part-ls} into \eqref{dou-sum-ls} yields
\BQN\label{ds-neg-ls}
\lim_{S\to\IF}\limsup_{u\to\IF}\frac1{v(u)\Psi(u)}\sum_{0\leq i<k\leq N_u-1}
q_{i,k}(u) = 0.
\EQN

Next, take $S=n$ for $n=2,3,\ldots$ and denote by $E_n$ the set of discontinuity points of
$\int_0^T\mathcal{B}_{\alpha,H(t)}^\eta (n_\eta,x)\td t$ on $(0,\mu_\eta([0,n_\eta]))$. For each $n\geq2$, $E_n$ has measure
$0$ since $\mathcal{B}_{\alpha,H(t)}^\eta(n_\eta,\cdot)$ is monotone in $x$ and uniformly bounded
for $t\in[0,T]$ by \nelem{upp-bou-BaeDx}.
Thus, in view of \eqref{low-upp-I1}, combing \eqref{wea-con-loc} with \eqref{ds-neg-ls} we get
\BQNY\nonumber
\limsup_{n\to\IF}\frac{\int_0^T \mathcal{B}_{\alpha,H(t)}^\eta (n_\eta,x)\td t}{n_\eta} &\leq& \liminf_{u\to\IF}\frac{\pk{ L_{\eta,u}^*[0,T] >x}}{v(u)\Psi(u)}\\\nonumber
&\leq& \limsup_{u\to\IF}\frac{\pk{ L_{\eta,u}^*[0,T] >x}}{v(u)\Psi(u)}\\\label{inequ-ls}
 &\leq& \liminf_{n\to\IF}\frac{\int_0^T \mathcal{B}_{\alpha,H(t)}^\eta (n_\eta,x)\td t}{n_\eta}
\EQNY
for any $x\in \{0\}\cup E^c$, where
 \BQN\label{def-Ec}
 E^c:=\R^+ \setminus \bigcup_{n=2}^{\IF} E_n.
  \EQN
Further, for all $t\in[0,T],x\geq0$ and any $S>1$
$$ \frac{\MB_{\alpha,H(t)}^{\eta}(\lfloor S\rfloor_\eta,x)}{\lceil S\rceil_\eta} \leq \frac{\MB_{\alpha,H(t)}^{\eta}(\Se,x)}{\Se} \leq\frac{\MB_{\alpha,H(t)}^{\eta}(\lceil S\rceil_\eta,x)}{\lfloor S\rfloor_\eta},$$
which implies that for any $x\in \{0\}\cup E^c$
 \BQN\nonumber
\lim_{u\to\IF} \frac{ \pk{ L_{\eta,u}^*[0,T] >x} }{v(u)\Psi(u)} &=& \lim_{n\to\IF} \frac{\int_0^T \mathcal{B}_{\alpha,H(t)}^\eta (n_\eta,x)\td t}{n_\eta}\\\nonumber
&=&\lim_{S\to\IF} \frac{\int_0^T \mathcal{B}_{\alpha,H(t)}^\eta (\Se,x)\td t}{\Se}\\\label{def-lim-B-ls}
&=&\lim_{S\to\IF} \frac{\int_0^T \mathcal{B}_{\alpha,H(t)}^\eta (S,x)\td t}{S}:=\MB_{\alpha}^{\eta,H}(x).
 \EQN

We determine $\MB_{\alpha}^{\eta,H}(x)$ by the right limit for each $x\in\bigcup_{n=2}^{\IF} E_n$. Hence, by  monotonicity, $\MB_{\alpha}^{\eta,H}(x)$
is well-defined for any $x\geq0$.
Let $x_0>0$ be any continuity point of $\MB_\alpha^{\eta,H}(\cdot)$.
Since $E^c$ is dense in $\R^+$ we can choose two sequences of points $\{y_n, z_n, n\in\N\}$ from $E^c$ such that
$y_n\nearrow x_0$ and $z_n\searrow x_0$. By the monotonicity again
\BQNY
\MB_\alpha^{\eta,H}(z_n)  =\liminf_{S\to\IF} \frac{\int_0^T \mathcal{B}_{\alpha,H(t)}^\eta (S,z_n)\td t}{S} &\leq& \liminf_{S\to\IF} \frac{\int_0^T \mathcal{B}_{\alpha,H(t)}^\eta (S,x_0)\td t}{S}\\
&\leq& \limsup_{S\to\IF} \frac{\int_0^T \mathcal{B}_{\alpha,H(t)}^\eta (S,x_0)\td t}{S}\\
 &\leq& \limsup_{S\to\IF} \frac{\int_0^T \mathcal{B}_{\alpha,H(t)}^\eta (S,y_n)\td t}{S}= \MB_\alpha^{\eta,H}(y_n),
\EQNY
and similarly
\BQNY
&&\MB_\alpha^{\eta,H}(z_n) =\liminf_{u\to\IF} \frac{ \pk{ L_{\eta,u}^*[0,T] >z_n} }{v(u)\Psi(u)}  \leq \liminf_{u\to\IF} \frac{ \pk{ L_{\eta,u}^*[0,T] >x_0} }{v(u)\Psi(u)}\\
  &&\qquad\quad\leq  \limsup_{u\to\IF} \frac{ \pk{ L_{\eta,u}^*[0,T] >x_0} }{v(u)\Psi(u)} \leq
 \limsup_{u\to\IF} \frac{ \pk{ L_{\eta,u}^*[0,T] >y_n} }{v(u)\Psi(u)}  = \MB_\alpha^{\eta,H}(y_n).
\EQNY
Letting  $n\to\IF$ in the above inequalities implies
that \eqref{def-lim-B-ls} holds also for any $x>0$ continuity point of $\MB_\alpha^{\eta,H}(\cdot)$.\\
 Next we show that $\MB_\alpha^{\eta,H}(\cdot)$ is finite and positive.  The finiteness follows from \nelem{upp-bou-BaeDx}
 in Appendix.
In order to prove positivity of $\MB_\alpha^{\eta,H}(\cdot)$, we note that by Bonferroni inequality
\BQNY\nonumber
\pk{ L_{\eta,u}^*[0,T] >x} &\geq& \pk{\bigcup_{k=0}^{\lfloor N_u/2\rfloor-1} \LT\{L_{\eta,u}^* \Delta_{2k} >x\RT\}}\\\label{positi-B-cons}
&\geq& \sum_{k=0}^{\lfloor N_u/2\rfloor-1} \pk { L_{\eta,u}^* \Delta_{2k} >x }  - \sum_{0\leq i<k\leq \lfloor N_u/2\rfloor-1}
q_{2i,2k}(u).
\EQNY
Let $x\in(0,\mu_\eta([0,S_\eta]))$ be a continuity point of $\int_0^T\mathcal{B}_{\alpha,H(t)}^\eta (\Se,x)\td t$, then by a similar argument as used in \eqref{wea-con-loc}
\BQNY
\lim_{u\to\IF}\frac1{v(u)\Psi(u)} \sum_{k=0}^{\lfloor N_u/2\rfloor-1} \pk { L_{\eta,u}^* \Delta_{2k} >x }
&=& \limit{u}  \frac 1 {2\Se }\int_0^T  F_u(t,x)  \mu_{2\Seu}(\td t)\notag \\\label{wea-con-ls-aux}
&=&		\frac{1}{2\Se} \int_0^T \mathcal{B}_{\alpha,H(t)}^\eta (\Se,x)\td t.
\EQNY
Further, as shown in \eqref{dou-sum-ls}-\eqref{negi-notfar-ls}
\BQNY
\limsup_{u\to\IF}\frac1{v(u)\Psi(u)} \sum_{0\leq i<k\leq \lfloor N_u/2\rfloor-1}
q_{2i,2k}(u) \leq
2\lceil (16\overline{h})^{1/\alpha}\rceil^2 H_{\alpha}^2([0,1])T \Se
\exp(-\frac{1}{16}\underline{h}(\Se)^{\alpha/2}).
\EQNY
Consequently,
\BQNY
&&\liminf_{u\to\IF} \frac{ \pk{ L_{\eta,u}^*[0,T] >x} }{v(u)\Psi(u)}\\
&&\qquad \geq \frac1{2\Se}\LT( \int_0^T \mathcal{B}_{\alpha,H(t)}^\eta (\Se,x)\td t
- 4\lceil (16\overline{h})^{1/\alpha}\rceil^2 H_{\alpha}^2([0,1])T \Se^2
\exp(-\frac{1}{16}\underline{h}(\Se)^{\alpha/2})\RT),
\EQNY
hence the proof follows. \QED

\def\LL{\rho}
\def\Qal{\mathcal{Q}_{\alpha}(x)}
\def\QaS{\mathcal{Q}_{\alpha}(S,x)}
\def\Qabl{\mathcal{Q}_{\alpha}^b(x)}
\def\LuL{L_{\eta,u}^* \Lambda_{u}}
\def\tNeu{\widetilde N_u^\eta}

\prooftheo{MainTheo}
First note that for any $\eta,x\geq0$
\BQNY\label{negi-proba}
\pk{ \LuL >x} \leq \pk{L_{\eta,u}^*[-T,T]>x} \leq  \pk{ \LuL >x} + \pk{\sup_{t\in[-T,T]\setminus\Lambda_u} X(t) >u },
\EQNY
where
$$ \delta(u)= \bigl(  \ln u/u\bigr) ^{2/\beta} \quad \textrm{and}\quad\Lambda_u=[-\delta(u),\delta(u)]. $$
By {\bf A0-A2}, for arbitrary $\vp_1>0$ there exist $\vp\in (0,T)$ such that
\BQN \label{holder-cond}
\E{ \zE{(} X(t)/\sigma(t)-X(s)/\sigma(s) \zE{)}^2} \leq 3\abs{t-s}^\alpha,\quad \forall\, s,t\in[-\vp,\vp],
\EQN
$$\sigma(t)\leq 1- (1-\vp_1)b\abs{t}^\beta,\quad \forall\, t\in[-\vp,\vp],$$
$$\sigma(t)\leq 1- (1-\vp_1)b\vp^\beta,\quad \forall\, t\in[-T,T]\setminus[-\vp,\vp].$$
Consequently, by Piterbarg  inequality (see e.g.,
\cite{Pit96}[Theorem 8.1]) 
for large enough $u$ and   some positive $C$
\begin{eqnarray}
 \pk{\sup_{t\in[-\vp,\vp]\setminus\Lambda_u} X(t) >u } \leq  2C\vp u^{2/\alpha-1}\exp\LT( -\frac{u^2}{2(1-(1-\vp_1)b\delta^\beta(u))^2} \RT). \
\label{b0}
\end{eqnarray}
By Borell-TIS inequality (see Theorem 2.1.1 in \cite{AdlerTaylor}) 
for some positive $C_1$
\begin{eqnarray}
\pk{\sup_{t\in[-T,T]\setminus[-\vp,\vp]} X(t) >u } \leq \exp\LT( -\frac{(u-C_1)^2}{2(1-(1-\vp_1)b\vp^\beta)^2} \RT).
\label{b1}
\end{eqnarray}
 Combing (\ref{b0}) with (\ref{b1}) we get
\BQN\label{neg-big-inter}
 \pk{\sup_{t\in[-T,T]\setminus\Lambda_u} X(t) >u } = o\LT( \Psi(u)\RT)
 \EQN
as $u\to\IF$.
Hence
$$ \pk{L_{\eta,u}^*[-T,T]>x} \sim \pk{\LuL>x}, \quad u\to \IF,$$
if the latter is asymptotically equivalent to $\Psi(u)$,
and thus we need to investigate the asymptotics of $\pk{\LuL>x}$.\\
\underline{Ad i)} We use the same notation as introduced in the proof of \netheo{theo-ls-tail}. Let
$$  \Delta_k=[k\Seu,(k+1)\Seu], \quad k=0,\pm1\ldot \pm N_u', $$
where $N_u'=\lfloor \delta(u)/\Seu\rfloor$. By Bonferroni inequality
\BQN\label{ineq-alp-s-b}
I_1'(u) \leq \pk{ \LuL  >x} \leq I_2'(u)
\EQN
holds for any $x\geq0$, where
\BQNY
&& I_1'(u) = \sum_{k=-N_u'}^{N_u'-1}\pk{ L_{\eta,u}^* \Delta_{k} >x }\ - \sum_{-N_u'\leq i<k\leq N_u'-1} q_{i,k}(u),\\
&& I_2'(u) = \sum_{k=-N_u'-1}^{N_u'}\pk{ L_{\eta,u}^* \Delta_{k} >x }\ + \sum_{-N_u'-1\leq i<k\leq N_u'} q_{i,k}(u).
\EQNY
Next, set
\BQNY
\xi_{u,k}(t)= \frac{ X\LT( k\Seu+t/v(u)\RT) }{ \sigma\LT( k\Seu+t/v(u)\RT)},\ t\in[0,\Se],
\EQNY
and
\BQNY
g_{k}(u)= \left\{\begin{array}{ll}u\big(1+(1-\vp_1)b\abs{k\Seu}^\beta \big),& k\in K_u, k\ge 0,\\
	u\big(1+(1-\vp_1)b\abs{(k+1)\Seu}^\beta \big),& k\in K_u, k<0,
\end{array}\right.
\EQNY
where $K_u= \{-N_u'-1 \ldot 0 \ldot  N_u'\}$.
It follows that  $g_{k}(u)$ converges as $u\to\IF$ to infinity uniformly for $ k\in K_u$. Moreover, assumptions {\bf C1-C3} in \netheo{the-weak-conv} are fulfilled by the family of Gaussian processes $\{\xi_{u,k}(t),t\in[0,\Se], k\in K_u\}$ given above. Specifically,
$h(t)=t^\alpha$ and $ \zeta(t)=B_{\alpha}(t)$ for $ t\in[0,\Se]$, and $\nu$ required in {\bf C3}
as shown by \eqref{holder-cond} is equal to $\alpha$. Therefore, by the uniform convergence as stated in Theorem \ref{the-weak-conv}, we have
\BQNY
\sum_{k=-N_u'-1}^{N_u'} \pk{ L_{\eta,u}^*\Delta_{k} >x } &\leq&
\sum_{k=-N_u'-1}^{N_u'} \pk{ \int_{[0,\Se]}\njk{g_{k}(u)(\xi_{u,k}(t)-g_{k}(u))}\mu_\eta(\td t) >x }\\
&\sim& \MB_{\alpha}^\eta(\Se,x) \sum_{k=-N_u'-1}^{N_u'} \Psi(g_{k}(u)),\quad u\to\IF
\EQNY
at $x=0$ and $x\in(0,\mu_\eta([0,S_\eta]))$ a continuity point of $\MB_{\alpha}^\eta(\Se,x)$,
where $\MB_{\alpha}^\eta(\Se,x)$=$\MB_{\alpha,1}^\eta(\Se,x)$ with the latter defined in \eqref{Baleta}.
Further, as $u\to \IF$,
\BQNY
\sum_{k=-N_u'-1}^{N_u'} \Psi(g_{k}(u))&\sim& \frac{2}{\sqrt{2\pi}u} \sum_{k=0}^{N_u'} \exp\LT( -\frac{u^2\LT( 1+(1-\vp_1)b\abs{k\Seu}^\beta\RT)^2}{2} \RT) \\
  &\sim& \frac{2\Psi(u)}{\Seu} \int_{0}^{\delta(u)} \exp\LT(-b(1-\vp_1)u^2t^{\beta}\RT)\td t \\
  & \sim& \frac{2(b(1-\vp_1))^{-1/\beta}}{\Se}\Gamma(1/\beta+1)u^{2/\alpha-2/\beta}\Psi(u)
\EQNY
and thus
\BQN\label{upp-onedex-ab}
 \limsup_{u\to\IF}\frac{\sum_{k=-N_u'-1}^{N_u'} \pk{ L_{\eta,u}^*(\Delta_{k}) >x }  }{2b^{-1/\beta}\Gamma(1/\beta+1)u^{2/\alpha-2/\beta}\Psi(u)}\leq (1-\vp_1)^{-1/\beta}  \frac{\MB_{\alpha}^\eta(\Se,x)}{\Se}
 \EQN
at $x=0$ and all continuity points $x\in(0,\mu_\eta([0,S_\eta]))$.
Moreover, as shown in \cite{Pit20} (see p. 22 therein)
\BQN\label{dou-sum-alb}
 \lim_{S\to\IF}\limsup_{u\to\IF}\frac{1}{u^{2/\alpha-2/\beta}\Psi(u)}\sum_{-N_u'-1\leq i<k\leq N_u'} q_{i,k}(u)=0,
\EQN
Consequently, substituting \eqref{upp-onedex-ab} and \eqref{dou-sum-alb} into \eqref{ineq-alp-s-b}, then taking $S=n$ for $n=2,3,\ldots$ yields
$$ \limsup_{u\to\IF}\frac{\pk{ \LuL  >x}}{2b^{-1/\beta}\Gamma(1/\beta+1)u^{2/\alpha-2/\beta}\Psi(u)}\leq (1-\vp_1)^{-1/\beta} \liminf_{n\to\IF} \frac{\MB_{\alpha}^\eta(n\eta,x)}{n\eta} $$
at any $x\in \{0\}\cup E^c$, with $E^c$ as defined in \eqref{def-Ec}. Here $E_n$ denotes the set of discontinuity points of $\MB_{\alpha}^\eta(n\eta,x)$ on $(0,\mu_\eta([0,n_\eta]))$.

Similarly, for any $x\in \{0\}\cup E^c$
$$ \liminf_{u\to\IF}\frac{\pk{ \LuL  >x}}{2b^{-1/\beta}\Gamma(1/\beta+1)u^{2/\alpha-2/\beta}\Psi(u)}\geq (1+\vp_1)^{-1/\beta} \limsup_{n\to\IF} \frac{\MB_{\alpha}^\eta(n\eta,x)}{n\eta}. $$
Since $\vp_1$ is arbitrary, then by the same argument as used in the proof of \netheo{theo-ls-tail} we have
$$ \lim_{u\to\IF}\frac{\pk{ \LuL  >x}}{2b^{-1/\beta}\Gamma(1/\beta+1)u^{2/\alpha-2/\beta}\Psi(u)} =
 \lim_{n\to\IF} \frac{\MB_{\alpha}^{\eta}(n\eta,x)}{n\eta}
=\lim_{S\to\IF} \frac{\MB_{\alpha}^{\eta}(S,x)}{S}:=\MB_{\alpha}^{\eta}(x)  $$
at $x=0$ and any $x>0$ a continuity point of $\MB_{\alpha}^{\eta}(\cdot)$. This together with \eqref{neg-big-inter} validates the claim \eqref{con-alp-s-bet-tail}.

\def\LuS{L_{u,S}}

\underline{Ad ii)} Set for large $S$
  \BQN\label{defi-del-S}
   \Delta_S=[-S/v(u),S/v(u)]
   \EQN
  and then for arbitrary $x\geq0$
\BQNY\label{ineq-alp-bet}
\pk{ L_{\eta,u}^* \Delta_S >x} \leq \pk{ \LuL >x} \leq \pk{ L_{\eta,u}^*\Delta_{S} >x} +
\pk{ \sup_{t\in\Lambda_u\setminus\Delta_{S}} X(t)>u}.
\EQNY
It follows from \nelem{lemm-soj-sta-sma}-ii) that
\BQNY
\lim_{u\to\IF} \frac{ \pk{ L_{\eta,u}^* \Delta_S >x} } {\Psi(u)} = \mathcal{P}_\alpha^{b,\eta}(S,x)
\EQNY
at $x=0$ and all continuity points $x\in(0,\mu_\eta([-S, S]))$ of $\mathcal{P}_\alpha^{b,\eta}(S,x)$ defined in \eqref{defi-B-b-eta-S}.
Further, as shown in \cite{Pit20} (see p. 22 therein),
\BQNY
\pk{ \sup_{t\in\Lambda_u\setminus\Delta_S} X(t)>u} = O\LT( e^{-cS^\alpha}\RT)\Psi(u)\oo, \quad
u\to\IF
\EQNY
holds for some $c>0$.
Then, with similar arguments as in the proof of \netheo{theo-ls-tail}, we obtain
$$\lim_{u\to\IF} \frac{\pk{ \LuL >x} } {\Psi(u)} =\lim_{S\to\IF} \mathcal{P}_\alpha^{b,\eta}(S,x) = \mathcal{P}_\alpha^{b,\eta}(x)\in(0,\IF) $$
at $x=0$ and any $x>0$ a continuity point of $\mathcal{P}_\alpha^{b,\eta}(\cdot)$.
The finiteness of $\mathcal{P}_\alpha^{b,\eta}(\cdot)$ follows from the fact that
$\mathcal{P}^{b,\eta}_\alpha(x)\leq \mathcal{P}_{\alpha}^b.$ Using further \eqref{neg-big-inter} establishes \eqref{LCohen}.

\underline{Ad iii)} For large $S$ define $\Delta_S$ as in \eqref{defi-del-S}.
Note that $v(u)=u^{2/\beta}$ since $\alpha>\beta$. For any $\vp>0$ and all large $u$, we have  $\delta(u)<\vp u^{-2/\alpha}$.
Hence for any $x\geq0$
\BQNY
\pk{ L_{\eta,u}^* \Delta_S >x} \leq \pk{ \LuL >x} \leq \pk{ L_{\eta,u}^*\Delta_{S} >x} +
\pk{ \sup_{t\in[-\vp u^{-2/\alpha}, \vp u^{-2/\alpha}]\setminus\Delta_S} X(t)>u}.
\EQNY
In view of \nelem{lemm-soj-sta-sma}-ii) we have
\BQNY
\lim_{u\to\IF} \frac{ \pk{ L_{\eta,u}^* \Delta_S >x} } {\Psi(u)} = \mathcal{T}_\beta^{b,\eta}(S,x)
\EQNY
at $x=0$ and all continuity points $x\in(0,\mu_\eta([-S, S]))$ of $\mathcal{T}_\beta^{b,\eta}(S,x)$ defined in \eqref{defi-T-b-eta-S}. Further, Lemma 5.1 in \cite{debicki2016parisian} implies
\BQNY
\pk{ \sup_{t\in[-\vp u^{-2/\alpha}, \vp u^{-2/\alpha}]\setminus\Delta_S} X(t)>u} &\leq&
\pk{ \sup_{t\in[-\vp u^{-2/\alpha}, \vp u^{-2/\alpha}]\setminus\Delta_S} \frac{X(t)}{\sigma(t)} > u\LT(1+(1-\vp)b\abs{\frac{S}{v(u)}}^\beta \RT) }\\
&\leq& \E{\sup_{s\in[-\vp,\vp]}e^{\sqrt{2}B_\alpha(s)-s^\alpha}}e^{-b(1-\vp)S^\beta}\Psi(u)\oo,
\quad u\to \IF.
\EQNY
Following the same argument as in case ii), we obtain
\BQNY
\lim_{u\to\IF} \frac{ \pk{ L_{\eta,u}^*[-T,T] >x} } {\Psi(u)} = \lim_{S\to\IF}\mathcal{T}_\beta^{b,\eta}(S,x):=\mathcal{T}_\beta^{b,\eta}(x)
\EQNY
at $x=0$ and all positive continuity points of $\mathcal{T}_\beta^{b,\eta}(\cdot)$, where for $\eta=0$, $\mathcal{T}^{b,0}_\beta (x)=e^{-b(\frac{x}{2})^\beta}$ if $x\geq0$ and for $\eta>0$,
$\mathcal{T}^{b,\eta}_\beta(x)=1$ if $x\in [0,\eta) $ and
$\mathcal{T}^{b,\eta}_\beta(x) =
e^{-b(k\eta)^\beta}$ if $x\in[(2k-1)\eta,(2k+1)\eta), k\in \N $.

 This completes the proof.  \QED

\section{Appendix}
\def\eZn{\eta\mathbb{Z}^n}
\def\medt{\mu_{\eta}(\td\vk{t})}
\def\xukt{\xi_{u,k}(\vk{t})}
\def\cukt{\chi_{u,k}(\vk{t})}
\def\txukt{\tilde\xi_{u,k}(\vk{t})}
\def\sukt{\sigma_{\xi_{u,k}}(\vk{t})}
\def\suks{\sigma_{\xi_{u,k}}(\vk{s})}
\def\suk0{\sigma_{\xi_{u,k}}(\vk{0})}
\def\siguk{\sigma_{\xi_{u,k}}}
\def\Ruk{R_{u,k}}
\def\ruk{\rho_{u,k}}
\def\vkt0{(\vk{t},\vk{0})}
\def\vks0{(\vk{s},\vk{0})}
\def\vkst{(\vk{s},\vk{t})}
\def\iukxz{\mathcal{I}_{u,k}(x;z)}

Let $K_u$ be an index function of $u$, $\vk{D}$ be a compact set in $\mathbb{R}^n$ and suppose without loss
 of generality that $\vk{0}\in\vk{D}$. Further, let
 $\{\xi_{u,k}(\vk{t}),\vk{t}\in\vk{D},k\in K_u\}$ be a family of centered Gaussian random fields with a.s. continuous sample paths and variance function
 $\sigma^2_{\xi_{u,k}}$. For $\vk{t}$ such that $\sigma^2_{\xi_{u,k}}(\vk{t})>0$ define the standardised  process
$$\tilde\xi_{u,k}(\vk{t}):= \frac \xukt { \sigma_{\xi_{u,k}}(\vk{t})},\quad \vk{t}\in\vk{D}.$$
Suppose that:
\begin{itemize}
 \item[{\bf C0}:]   $\{g_{k}(u),k\in K_u\}$ is a sequence of deterministic functions of $u$ satisfying
     \BQNY
     \lim_{u\to\IF}\inf_{k\in K_u}g_{k}(u)=\IF.
     \EQNY
\item[{\bf C1}:] $\sigma_{\xi_{u,k}}(\vk{0})=1$ for all large $u$ and any $k\in K_u$, and there exists some bounded continuous function $h$ on $\vk{D}$ such that
\BQNY\label{assump-cova-field}
\lim_{u\to\IF}\sup_{\vk{t}\in \vk{D},k\in K_u}\abs{g_{k}^2(u)\LT( 1- \E{ \xi_{u,k}(\vk{t})\xi_{u,k}(\vk{0}) } \RT) - h(\vk{t}) } =0.
\EQNY	
\item[{\bf C2}:] There exists a centered Gaussian random field $\zeta(\vk{t}),\vk{t}\in\mathbb{R}^n$ with a.s. continuous trajectories such that for any $\vk{s},\vk{t}\in\vk{D}$
\BQNY\label{assump-vari-field}
  \lim_{u\to\IF}\sup_{k\in K_u}\abs{g_{k}^2(u)\big(Var(\tilde\xi_{u,k}(\vk{t})-\tilde\xi_{u,k}(\vk{s}))\big) - 2Var(\zeta(\vk{t})-\zeta(\vk{s})) } =0.
\EQNY
\item[{\bf C3}:] There exist positive constants $C, \nu, u_0$ such that
\BQNY\label{assump-holder-field}
  \sup_{k\in K_u}  g_{k}^2(u)\E{(\xi_{u,k}(\vk{t})-\xi_{u,k}(\vk{s})\Ehb{)}^2} \leq C \norm{\vk{s}-\vk{t}}^\nu
  \EQNY
  holds for all $\vk{s},\vk{t}\in\vk{D}, u\geq u_0$, where $\norm{\vk{t}}^v=\sum_{i=1}^n\abs{t_i}^v$.
\end{itemize}
We present below an extension of Theorem 2.1 in \cite{KEP2016}. Hereafter, $C_i, i\in \N$ are  positive constants which might be different from line to line.
We recall that
$\medt/\eta^n$ denotes the counting measure on $\eZn, \eta>0$ and
$\mu_0$ is the Lebesgue measure on $\R^n$.

\BT\label{the-weak-conv}
Let $h,g_{k}(u),\xi_{u,k}(\vk{t}),\vk{t}\in\vk{D},k\in K_u$ and $\zeta$ be such that  {\bf C0-C3} hold.
Then, for $\eta\ge 0$
\BQN\label{con-uni-con}
\lim_{u\to\IF}\sup_{k\in K_u} \abs{ \frac{\pk{ \int_{\vk{D}} \mathbb{I}_0\LT(g_{k}(u)(\xi_{u,k}(\vk{t})-g_{k}(u))\RT) \medt >x }  } {\Psi(g_{k}(u))} -\MB^{h,\eta}_{\zeta}(\vk{D},x) } = 0
\EQN
at $x=0$ and all $x\in(0,\mu_\eta(\vk{D}))$
continuity points of $\MB^{h,\eta}_{\zeta}(\vk{D},x)$, where
\BQNY
\MB^{h,\eta}_{\zeta}(\vk{D},x)= \int_\R  \pk{ \int_{\vk{D}} \mathbb{I}_0\big(\sqrt{2}\zeta(\vk{t})-h(\vk{t}) + z \big) \medt > x } e^{-z} \td z.
\EQNY
\ET

\prooftheo{the-weak-conv}
Suppose that {\bf C0}-{\bf C3} are satisfied. We begin from the observation that
  \BQN\label{assump-holder-field-corr}
 \limsup_{u\to\IF}\sup_{k\in K_u}  g_{k}^2(u)\E{(\tilde\xi_{u,k}(\vk{t})-\tilde\xi_{u,k}(\vk{s}))^2} \leq C_1 \norm{\vk{s}-\vk{t}}^\nu,\quad \forall\, \vk{s},\vk{t}\in\vk{D},
  \EQN
where $C_1,\nu$ are positive constants.
Indeed, note that
$$1-\sigma^2_{\xi{u,k}}(\vk{t})=2\LT(1- \E{ \xi_{u,k}(\vk{t})\xi_{u,k}(\vk{0})} \RT) - \E{(\xukt-\xi_{u,k}(\vk{0}))^2},$$
which together with {\bf C1} and {\bf C3} implies
\BQN\label{varia-uni-conv}
\lim_{u\to\IF}\sup_{\vk{t}\in \vk{D},k\in K_u}\abs{\sigma^2_{\xi_{u,k}}(\vk{t})-1}=0.
\EQN
Consequently, for sufficiently large $u$
    \BQNY
     g_{k}^2(u)\E{(\tilde\xi_{u,k}(\vk{t})-\tilde\xi_{u,k}(\vk{s}))^2} &=&g_{k}^2(u)\frac{2\sukt\suks-2\E{ \xi_{u,k}(\vk{t})\xi_{u,k}(\vk{s}) } }{\sukt\suks}\\
     &\leq& g_{k}^2(u)\frac{\E{(\xi_{u,k}(\vk{t})-\xi_{u,k}(\vk{s}))^2}} {\inf_{\vk{t}\in\vk{D}}\sigma^2_{u,k}(\vk{t})} \leq 2C\norm{\vk{s}-\vk{t}}^\nu,\quad \forall\, k\in K_u, \vk{s},\vk{t}\in\vk{D}.
    \EQNY
Next, for notational simplicity denote by $R_{u,k}$ and  $\rho_{u,k}$ the covariance and the correlation function of $\xi_{u,k}$. Further set
$$ \chi_{u,k}(\vk{t}) := g_{k}(u)(\tilde\xi_{u,k}(\vk{t})-\rho_{u,k}(\vk{t},\vk{0})\tilde\xi_{u,k}(\vk{0})),\quad \vk{t}\in\vk{D}
$$
and
$$f_{u,k}(\vk{t},z):=z R_{u,k}(\vk{t},\vk{0})-g_{k}^2(u)\LT( 1-R_{u,k}(\vk{t},\vk{0}) \RT)  ,\ \vk{t}\in \vk{D}, z\in\R.$$
Conditioning on $\xi_{u,k}(\vk{0})$ and using that $\xi_{u,k}(\vk{0})$ and $\xi_{u,k}(\vk{t})-\Ruk(\vk{t},\vk{0})\xi_{u,k}(\vk{0})$ are mutually independent for large $u$, we obtain
\BQNY
\lefteqn{\pk{ \int_{\vk{D}} \njk{g_{k}(u)(\xi_{u,k}(\vk{t})-g_{k}(u))} \medt >x}}\\
&=&\frac{e^{-g_{k}^2(u)/2}} { \sqrt{2\pi}g_{k}(u)}
\int_{\R} \exp\LT(-z-\frac{z^2}{2g_{k}^2(u)}\RT)\\
 && \times 
 \pk{ \int_{\vk{D}} \njk{g_{k}(u)(\xi_{u,k}(\vk{t})-g_{k}(u))} \medt >x | \xi_{u,k}(\vk{0})=g_{k}(u)+zg_{k}^{-1}(u)} \td z\\
&=&\frac{e^{-g_{k}^2(u)/2}} { \sqrt{2\pi}g_{k}(u)}
\int_{\R} \exp\LT(-z-\frac{z^2}{2g_{k}^2(u)}\RT) \pk{ \int_{\vk{D}} \njk{\sukt\chi_{u,k}(\vk{t}) + f_{u,k}(\vk{t},z) } \medt >x } \td z.
\EQNY
Let
\BQNY
\mathcal{I}_{u,k}(x;z):=\pk{ \int_{\vk{D}} \njk{\sukt\chi_{u,k}(\vk{t}) + f_{u,k}(\vk{t},z) } \medt >x }.
\EQNY
Consequently, in order to show the claim it suffices to prove that
\BQN\label{aux-uni-conv}
\lim_{u\to\IF}\sup_{k\in K_u}\abs{ \int_{\R} \exp\LT(-z-\frac{z^2}{2g_{k}^2(u)}\RT)  \mathcal{I}_{u,k}(x;z) \td z - \MB^{h,\eta}_{\zeta}(\vk{D},x)} =0
\EQN
at $x=0$ and all $x\in(0,\mu_\eta(\vk{D}))$ positive continuity points of $\MB^{h,\eta}_{\zeta}(\vk{D},x)$. Since for all $x\geq0$ and any large $M$
$$ \sup_{k\in K_u}e^{-z}\iukxz \leq e^{-z},\quad z\geq-M$$
and by Piterbarg inequality for all large $u$ and $M$
\BQN\nonumber
\sup_{k\in K_u}e^{-z}\iukxz &\leq&
\sup_{k\in K_u}\pk{\sup_{\vk{t}\in\vk{D}}\{\sukt\chi_{u,k}(\vk{t}) + f_{u,k}(\vk{t},z)\}>0} e^{-z}\\ \nonumber
 &\leq& \sup_{k\in K_u}\pk{\sup_{\vk{t}\in\vk{D}} \chi_{u,k}(\vk{t})>C_2\abs{z}-C_3} e^{-z}\\ \label{for-domi-Iuk}
 &\leq& C_4 \abs{z}^{2n/\nu-1}e^{-C_5z^2-C_6z},\quad z<-M,
\EQN
then by the dominated convergence theorem and assumption {\bf C0}
\BQNY
\lefteqn{\sup_{k\in K_u}\abs{\int_{\R}  \exp\LT(-z-\frac{z^2}{2g_{k}^2(u)} \RT)  \iukxz \td z -\int_{\R}e^{-z}\iukxz \td z }}\\
&&\leq \int_{\R} \sup_{k\in K_u} \LT(e^{-z} \iukxz\RT) \abs{1- e^{-z^2/(2g_{k}^2(u))}  } \td z \rightarrow 0, \quad u\to\IF.
\EQNY
Therefore, in order to prove the convergence in \eqref{aux-uni-conv} it suffices to show that
\BQN\label{aux-2-uni-conv}
\lim_{u\to\IF}\sup_{k\in K_u}\abs{\int_{\R} e^{-z} \mathcal{I}_{u,k}(x;z) \td z - \MB^{h,\eta}_{\zeta}(\vk{D},x) } =0
\EQN
at $x=0$ and all continuity points $x\in(0,\mu_\eta(\vk{D}))$.

 Let  $C(\vk{D})$ denote the Banach space of all continuous functions on $\vk{D}$ equipped with sup-norm.  For any $\vk{s},\vk{t}\in\vk{D}$, from {\bf C2} and \eqref{assump-holder-field-corr} we have
\BQNY Var(\chi_{u,k}(\vk{t})-\chi_{u,k}(\vk{s})) &=& g_{k}^2(u)\LT( \E{(\tilde\xi_{u,k}(\vk{t})-\tilde\xi_{u,k}(\vk{s})\Ehb{)^2}} - \LT( \rho_{\xi_{u,k}}(\vk{t},\vk{0})-\rho_{\xi_{u,k}}(\vk{s},\vk{0}) \RT)^2 \RT) \\&\to & 2Var(\zeta(\vk{t})-\zeta(\vk{s}))
\EQNY
uniformly with respect to $k\in K_u$ as $u\to\IF$.
 Hence, the finite-dimensional distributions of $\chi_{u,k}$ converge to that
 of $\sqrt{2}\zeta(\vk{t}),\vk{t}\in\vk{D}$ uniformly with respect to $k\in K_u$. In view of \eqref{assump-holder-field-corr},
 we know that
 the measures on $C(\vk{D})$ induced by $\{\chi_{u,k}(\vk{t}),\vk{t}\in\vk{D},k\in K_u\}$ are uniformly tight for large $u$, and by \eqref{varia-uni-conv} $\sukt$ converges to $1$ uniformly for $\vk{t}\in \vk{D}$ and $k\in K_u$ as $u\to \IF$.
 Therefore,  $\{\sukt\chi_{u,k}(\vk{t}),\vk{t}\in\vk{D}\}$
 converge weakly to $\{\sqrt{2}\zeta(\vk{t}),\vk{t}\in\vk{D}\}$ as $u\to\IF$ uniformly with respect to $k\in K_u$. Further, by {\bf C0}-{\bf C1} for each $z\in\Z$
 $$ \limit{u} \sup_{k\in K_u, \vk{t}\in\vk{D}}\abs{f_{u,k}(\vk{t},z)-z +h(\vk{t})}
 =0 $$
 implying that for each $z\in\Z$, the probability measures on $C(\vk{D})$ induced by
 $\{\chi_{u,k}^f(\vk{t},z),\vk{t}\in\vk{D}\}$,
where
 $$\chi_{u,k}^f(\vk{t},z):=\sukt\chi_{u,k}(\vk{t}) + f_{u,k}(\vk{t},z)\quad \textrm{and}\quad \zeta_h(\vk{t}):=\sqrt{2}\zeta(\vk{t})-h(\vk{t}),$$
converge weakly, as $u\to\IF$, to that induced by $\{\zeta_h(\vk{t}) + z,\vk{t}\in\vk{D}\}$ uniformly with respect to $k\in K_u$,

{Consequently, for any $\eta>0,z\in\Z$
	\BQN\label{wea-conv-cont-aux}
	\lim_{u\to\IF}\sup_{k\in K_u}\abs{ \mathcal{I}_{u,k}(x,z)  - \mathcal{I}(x;z)} =0
	\EQN
	holds at all continuity points $x\in(0,\mu_\eta(\vk{D}))$ (depending on $z$) of $\mathcal{I}(x;z)$  defined  by
	$$\mathcal{I}(x;z) := \pk{ \int_{\vk{D}} \njk{ \zeta_h(\vk{t}) + z } \medt > x }. $$
	For $\eta=0$, by \cite{berman1973excursions}[Lemma 4.2] the set of
discontinuity points  of
	$$\int_{\vk{D}} \mathbb{I}_0(f(\vk{t})) d\vk{t},\quad f\in C(\vk{D})$$
	is of measure $0$ under the probability measure induced by $\{\zeta_h(\vk{t}) + z,\vk{t}\in\vk{D}\}$. Consequently, by the continuous mapping theorem we also have \eqref{wea-conv-cont-aux}.}
Next, we borrow an argument from \cite{Berman92}[Theorem 1.3.1] to verify \eqref{aux-2-uni-conv} for all positive continuity points.
 Let $x\in(0,\mu_\eta(\vk{D}))$ be such a continuity point, i.e.,
$$ \lim_{\vp\to0}\int_{\R} \LT(\mathcal{I}(x_0+\vp;z) -\mathcal{I}(x_0-\vp;z)\RT)e^{-z} \td z=0.$$
Since for large $M$ and all $x\geq0$ by Borell-TIS inequality
\BQN\label{for-domi-IF}
e^{-z}\mathcal{I}(x;z) \leq C_7 e^{-C_8z^2-C_9z},\quad z<-M
\EQN
it follows from the dominated convergence theorem that
 $$ \int_{\R} \LT(\mathcal{I}(x_0+;z) -\mathcal{I}(x_0-;z)\RT) e^{-z} \td z=0,$$
and thus by the monotonicity of $\mathcal{I}(x;z)$ in $x$ for each fixed $z$, $x_0$ is a continuity point of $\mathcal{I}(x;z)$ for almost all $z\in\R$.
Hence by \eqref{wea-conv-cont-aux} for almost all $z\in\R$
\BQNY\label{exp-almost-z}
 \lim_{u\to\IF}\sup_{k\in K_u}\abs{ \mathcal{I}_{u,k}(x_0,z)  - \mathcal{I}(x_0;z)} =0.
\EQNY
As shown in \eqref{for-domi-Iuk} and \eqref{for-domi-IF}
it follows from the dominated convergence that
\BQNY
\sup_{k\in K_u}\abs{\int_{\R} e^{-z} \mathcal{I}_{u,k}(x_0;z) \td z - \int_{\R}e^{-z} \mathcal{I}(x_0;z)\td z }
\leq\int_{\R} \sup_{k\in K_u}\abs{\mathcal{I}_{u,k}(x_0;z) -  \mathcal{I}(x_0;z)} e^{-z} \td z  \rightarrow 0,
\EQNY
as $u\to\IF$, establishing the proof for all continuity points $x\in(0,\mu_\eta(\vk{D}))$.\\
The case $x=0,\eta=0$ is shown in \cite{KEP2016}.
Since the case $x=0,\eta>0$ can be established
by arguments similar to the presented above, we omit the details.
This completes the proof.   \QED

Let for any $\eta\geq0,\,S>0$, $x\in[0,\mu_\eta([-S,S]))$
\BQN\label{defi-B-b-eta-S}
\mathcal{P}_\alpha^{b,\eta}(S,x) := \int_{\R} \pk{\int_{-S}^S
	\njk{ W_{\alpha}(s)-b\abs{s}^\alpha+ z } \mu_{\eta}(\td s) >x } e^{-z}\td z
\EQN
and
\BQN\label{defi-T-b-eta-S}
\mathcal{T}_\beta^{b,\eta}(S,x) \coloneqq \int_{0}^\IF \pk{\int_{-S}^S
	\njk{ -b\abs{s}^\beta+ z } \mu_{\eta}(\td s) >x } e^{-z}\td z.
\EQN

\BEL\label{lemm-soj-sta-sma}
i) Let $X$ be as in \netheo{theo-ls-tail} and let $v(u)$ be as in \eqref{defini-vu-loc}. For any $\eta\geq0,$ $S>\eta$ and $t_u,u>0$ such that $\limit{u} t_u=t_0\in [0,T]$, we have
\BQNY\label{weak-conv-loc-sta}
\lim_{u\to\IF} \Psi^{-1}(u) \pk{ v(u)\int_0^{S/v(u)}  \nj{X(t_u+s)} \mu_{\eta_u}(\td s)>x}
=  \mathcal{B}_{\alpha,H(t_0)}^\eta (S,x)
\EQNY
at $x=0$ and any $x\in(0,\mu_\eta([0, S]))$ continuity point of $\mathcal{B}_{\alpha,H(t_0)}^\eta (S,x)$. \\
ii) Let $X$ be as in \netheo{MainTheo} and $v(u)$ be defined in \eqref{defi-v}. Then for any $\eta\geq0,\,S>0$
\BQN\label{def-Babe-x}
\pk{ L^\ast_{\eta,u}[-S/v(u),S/v(u)] >x} \sim
\Psi(u)\times
\left\{
\begin{array}{ll}
\mathcal{P}_\alpha^{b,\eta}(S,x) & \mbox{if $\alpha=\beta$}\\[0.5cm]
\mathcal{T}_\beta^{b,\eta}(S,x) & \mbox{if $\alpha>\beta$} ,
\end{array}
\right.
\EQN
as $u\to \IF$,
for $x=0$ and $x\in(0,\mu_\eta([-S, S]))$ a continuity point of $\mathcal{P}_\alpha^{b,\eta}(S,x)$ or
$\mathcal{T}_\beta^{b,\eta}(S,x)$ respectively.

\EEL

\prooflem{lemm-soj-sta-sma}
i) For any $x\geq0$
\BQNY \pk{ v(u)\int_0^{S/v(u)}  \nj{X(t_u+s)} \mu_{\eta_u}(\td s)>x}
= \pk{ \int_0^S \njk{u(X(t_u+t/v(u))-u)} \mu_{\eta}(dt)>x }.
\EQNY
Set $\vk{D}=[0,S\Ehb{]}, g_{k}(u)=u, K_u=1$ and $\xi_{u,k}(t)=X(t_u+t/v(u))$. By the Uniform Convergence Theorem and Potter's Theorem (see e.g., \cite{Bingham1989}[Theorem 1.5.2 and Theorem 1.5.3 ]) it follows that  $\xi_{u,k}$ satisfies  the assumptions {\bf C1-C3} with
$$h(t)=H(t_0)\abs{t}^{\alpha},\ \zeta({t})=\sqrt{H(t_0)}B_{\alpha}(t)\
\textrm{for}\ t\in [0,S],\ C=C_\alpha\ \textrm{and}\ \nu=\alpha/2.$$
Hence the claim  follows by \netheo{the-weak-conv} with
$  \MB^{h,\eta}_{\zeta}(\vk{D},x)=\mathcal{B}_{\alpha,H(t_0)}^\eta (S,x)$ and the claim in ii) follows with similar arguments.  \QED

\BEL\label{upp-bou-BaeDx}
If $h$ and $\zeta$ given in {\bf C1-C2} satisfy  $h(\vk{t})=\Var\zeta(\vk{t})$
and $\zeta(\vk{t})=\sum_{i=1}^n \sqrt{\lambda_i}B_\alpha^{(i)}(t_i)$ for some
 positive constants $\lambda_1 \ldot \lambda_n$, where $B_\alpha^{(i)}$'s are independent fBm's with Hurst index $\alpha/2$,
then for any $x,\eta\geq0$ and $\vk{D}=\prod_{i=1}^n[0,T_i]$ we have
$$\mathcal{B}_{\zeta}^{h,\eta} (\vk{D},x)\leq \prod_{i=1}^n\lceil\lambda_i^{1/\alpha}\rceil\lceil T_i\rceil H_\alpha^n([0,1]).$$
\EEL

\prooflem{upp-bou-BaeDx}
Let $\xi$ be a mean zero homogeneous Gaussian field with covariance function  $c(\vk{s}+\vk{t}, \vk{s})=r(\vk{t})=\exp(-\sum_{i=1}^n\lambda_i\abs{t_i}^\alpha)$. Taking $ K_u=1,g_{k}(u)=u$ and $\xi_{u,k}(\vk{t})=\xi(u^{-2/\alpha}\vk{t})$  \netheo{the-weak-conv} yields  for $\eta\ge 0$, $x=0$ and $x\in(0,\mu_\eta(\vk{D}))$ a continuity point of the constant below
$$ \limit{u}\Psi^{-1}(u)\pk{ \int_{\vk{D}} \mathbb{I}_0\LT(u(\xi(u^{-2/\alpha}\vk{t})-u)\RT) \medt >x } = \MB^{h,\eta}_{\zeta}(\vk{D},x).$$
By the homogeneity of $\xi$, we have further
\BQNY
\lefteqn{\Psi^{-1}(u)\pk{ \int_{\prod_{i=1}^n[0,T_i]} \mathbb{I}_0\LT(u(\xi(u^{-2/\alpha}\vk{t})-u)\RT) \medt >x }}\\
 &\leq&
\Psi^{-1}(u) \pk{\sup_{\vk{t}\in\prod_{i=1}^n[0,T_i]}\xi(u^{-2/\alpha}\vk{t})>u} \\
&\leq& \Psi^{-1}(u) \prod_{i=1}^n\lceil T_i\rceil \pk{\sup_{\vk{t}\in[0,1]^n}\xi(u^{-2/\alpha}\vk{t})>u}\\
&\to& \prod_{i=1}^n\lceil T_i\rceil \E{\sup_{\vk{t}\in[0,1]^n} e^{\sqrt{2}\zeta(\vk{t})-h(\vk{t}) }} \\
&=&  \prod_{i=1}^n\lceil T_i\rceil \E{\sup_{t_i\in[0,\lambda_i^{1/\alpha}]} e^{  \sqrt{2}B_\alpha^{(i)}(t_i)-\abs{t_i}^\alpha }}\\
&\leq&  \prod_{i=1}^n\lceil\lambda_i^{1/\alpha}\rceil\lceil T_i\rceil H_\alpha^n([0,1])
\EQNY
as $u\to\infty$, where the last inequality follows from the fact $H_\alpha([0,T])\leq\lceil T\rceil H_\alpha([0,1])$. \QED

\BEL\label{lemm-negli-ds}
If $X$ is a centered Gaussian process fulfilling \eqref{locall-station}, $v(u)$ and $\vp$ are defined in \eqref{defini-vu-loc} and in \eqref{defi-vp-con1}, respectively, then for $0\leq S_1<S_2<T_1<T_2<\infty, T_1-S_2\geq 1$ and $u$ large enough such that
\BQN\label{defi-vp-con3}
(T_2-S_1)/v(u)\leq\vp,
\EQN
we have
\BQNY
\pk{\sup_{t\in[S_1,S_2]/v(u)}X(t)>u,\sup_{t\in[T_1,T_2]/v(u)}X(t)>u}\leq C(\alpha,S_1,S_2,T_1,T_2)\Psi(u),
\EQNY
where
$$C(\alpha,S_1,S_2,T_1,T_2)=2\lceil (16\overline{h})^{1/\alpha}\rceil^2 \lceil S_2-S_1\rceil\lceil T_2-T_1\rceil H_{\alpha}^2([0,1]) \exp(-\frac{1}{16}\underline{h}\abs{T_1-S_2}^{\alpha/2}),$$
with $\underline{h}=\inf_{t\in[0,T]}H(t)>0$ and $\overline{h}=\sup_{t\in[0,T]}H(t)<\infty.$
\EEL

\prooflem{lemm-negli-ds} We borrow some arguments from the proof of \cite{mi:17}[Lemma 5]. Define next
\begin{align*}
A_u&=[S_1,S_2]/v(u),\qquad\qquad\qquad\, B_u=[T_1,T_2]/v(u),\\
Y(s,t)&=X(s)+X(t),\qquad\qquad\sigma^2(s,t)=\Var (Y(s,t)).
\end{align*}
By \eqref{locall-station}, for sufficiently close $s,t\in[0,T]$
$$\frac{1}{2}\underline{h}K(\abs{s-t})\leq 1-\rho(s,t)\leq2\overline{h}K(\abs{s-t}).$$
Consequently, for sufficiently large $u$ such that
\eqref{defi-vp-con3} holds,
by \eqref{defi-vp-con1}
$$\inf_{(s,t)\in A_u\times B_u}\sigma^2(s,t)\geq 4-4\overline{h}\sup_{(s,t)\in A_u\times B_u}K(\abs{s-t})>2$$
and
$$\sup_{(s,t)\in A_u\times B_u}\sigma^2(s,t)\leq 4-\underline{h}\inf_{(s,t)\in A_u\times B_u}K(\abs{s-t}),$$
implying that
\BQN\nonumber
\pk{\sup_{t\in A_u}X(t)>u,\ \sup_{t\in B_u}X(t)>u}&\leq& \pk{ \sup_{(s,t)\in A_u\times B_u)} Y(s,t)>2u}\\\label{upp-dou-ineq1}
&\leq& \pk{\sup_{(s,t)\in A_u\times B_u}Y^*(s,t)>u^*},
\EQN
where $Y^*(s,t)=Y(s,t)/\sigma(s,t)$ and
\BQNY\label{defi-u-star}
u^*=\frac{2u}{\sqrt{4-\underline{h}\inf_{(s,t)\in A_u \times B_u}K(\abs{s-t})}}.
\EQNY
As in \cite{mi:17}, we have $(s_1,t_1)\in A_u\times B_u, (s_2,t_2)\in A_u\times B_u$ and $u$ sufficiently large
\BQNY
\Cov(Y^*(s_1,t_1),Y^*(s_2,t_2))&\geq& 1-8\overline{h}K(\abs{s_2-s_1})-8\overline{h}K(\abs{t_2-t_1}).
\EQNY
Let
$Z(s,t):=\frac{1}{\sqrt{2}}(\vartheta_1(s)+\vartheta_2(t)),$
where $\vartheta_i,i=1,2$ are mutually independent copies of a
mean zero stationary Gaussian process $\vartheta$ with unit variance and covariance function satisfying
$1-\Cov(\vartheta(s),\vartheta(t))\sim 32\overline h K(\abs{s-t}), \abs{t-s}\to0.$
As mentioned in the proof of Theorem \ref{theo-ls-tail}, the existence of such a Gaussian process is guaranteed by the Assertion in \cite{VP1}[p.265].
Hence by Slepian inequality, for sufficiently large $u$ we have
\BQN\nonumber
\pk{\sup_{(s,t)\in A_u\times B_u}Y^*(s,t)>u^*} &\leq& \pk{ \sup_{(s,t)\in A_u\times B_u} Z(s,t)>u^*}\\\label{upp-dou-ineq2}
&\leq& \pk{ \sup_{(s,t)\in(A'_u\times B'_u)} Z(s,t)>u^*},
\EQN
where the last inequality follows from stationarity with $A'_u=[0,(S_2-S_1)]/v(u)$ and $B'_u=[0,(T_2-T_1)]/v(u)$.\\
Next, set $\vk{D}=[0,S_2-S_1]\times[0,T_2-T_1], g_{k}(u)=u^*$ and $\xi_{u,k}(s,t)=Z(s/v(u),t/v(u))$.
It is straightforward to check that assumptions {\bf C0-C2} are fulfilled with
$$h(s,t)=16\overline{h}(\abs{s}^{\alpha}+\abs{t}^{\alpha})), \quad \zeta(s,t)=4\sqrt{\overline{h}}(B_{\alpha}^{1}(s)+B_{\alpha}^{2}(t)),$$
 where $B_{\alpha}^{(i)},i=1,2$ are two independent fBm's with Hurst index $\alpha/2$. Further, by Potter's Theorem, assumption {\bf C3} holds for $v=\alpha/2$ and some constant $C$ depending on the sides length of $\vk{D}$.
 Thus, by \netheo{the-weak-conv} and \nelem{upp-bou-BaeDx}, for u sufficiently large
\BQN\label{upp-dou-ineq3}
\pk{ \sup_{(s,t)\in A'_u\times B'_u} Z(s,t)>u^*}  \leq 2\lceil (16\overline{h})^{1/\alpha}\rceil^2 \lceil S_2-S_1\rceil\lceil T_2-T_1\rceil H_{\alpha}^2([0,1])\Psi(u^*).
\EQN
Moreover, since $T_1-S_2\geq1$, then by Potter's Theorem again, we have for
sufficiently large $u$, that
$$u^2\inf_{(s,t)\in A_u\times B_u}K(\abs{s-t}))\geq \frac{1}{2}\abs{T_1-S_2}^{\alpha/2},$$
which implies that
\BQNY
(u^*)^2=\frac{4u^2}{4-\underline{h}\inf_{(s,t)\in A_u \times B_u}K(\abs{s-t})}\geq u^2(1+\frac{1}{4}\underline{h}\inf_{(s,t)\in A_u\times B_u}K(\abs{s-t}))
\geq u^2+\frac{1}{8}\underline{h}\abs{T_1-S_2}^{\alpha/2}.
\EQNY
Consequently, the claim follows by
\eqref{upp-dou-ineq1}-\eqref{upp-dou-ineq3} and the fact that
$\sqrt{2\pi}\Psi(u)\leq u^{-1}e^{-\frac{1}{2}u^2}$ for $u>0$.  \QED

\BEL \label{Tilt} Let $W$ be an $N(0,1)$ random variable independent of $\mathcal{Z}$ which is exponentially distributed with parameter 1.  For any $c>0$ we have
\BQNY
\pk{ cW- c^2/2 + Z> 0}= 2 \pk{W > c/2}.
\EQNY
\EEL

\prooflem{Tilt} Since $Z>0$ almost surely, then
$$\pk{ cW- c^2/2 + Z> 0, cW- c^2/2 \ge 0 } = \pk{ W> c/2}.$$
Let the random variable $V$ be such that
$$\pk{V \le x}= \E{ e^{c W- c^2/2} \mathbb{I}( cW- c^2/2 \le x)}, \quad  x\inr.$$
It is well-known, see e.g., \cite{Htilt}[Lemma 7.1] that $V$ has an $N(c^2/2,c^2)$ distribution. Hence by the independence of
$Z$ and $W$
\BQNY
\pk{ cW- c^2/2 + Z> 0, cW- c^2/2 \le 0 }
&=&\E{  e^{cW- c^2/2} \mathbb{I}( cW- c^2/2\le 0)}\\
&=& \pk{V \le 0}= \pk{ W \le -   c/2}
\EQNY
establishing the proof. \QED

{\bf Acknowledgement}:
K.D. was partially supported by    NCN Grant No 2015/17/B/ST1/01102 (2016-2019). X.P. thanks the Fundamental Research Funds for the Central Universities (ZYGX2015J102) and National Natural Science Foundation of China (71501025,11701070) for partial financial support.
Financial support from the Swiss National Science Foundation Grant 200021-175752/1 is also kindly acknowledged.

\bibliographystyle{ieeetr}
\bibliography{BermanC}
\end{document}